\documentclass{ws-procs975x65}
\usepackage{subfigure,tikz,float} 
\newcommand{\Dom}{\mathcal{D}}
\newcommand{\elliptic}{\vartheta}
\renewcommand{\L}{\Lambda}
\newcommand{\drm}{\mathrm{d}}

\newcommand{\EE}{\mathbb{E}}
\newcommand{\NN}{\mathbb{N}}
\newcommand{\PP}{\mathbb{P}}
\newcommand{\RR}{\mathbb{R}}
\newcommand{\ZZ}{\mathbb{Z}}
\newcommand{\cD}{\mathcal{D}}
\newcommand{\cL}{\mathcal{L}}
\newcommand{\ra}{\rangle}
\newcommand{\la}{\langle}
\newcommand{\tL}{Q}
\newcommand{\Tr}{{\mathop{\mathrm{Tr} \,}}}
\DeclareMathOperator{\supp}{\operatorname{supp}}

\DeclareMathOperator{\dist}{\operatorname{dist}}
\begin{document}
\title{Equidistribution estimates for eigenfunctions and eigenvalue bounds for random operators}
\author{D. Borisov}
\address{Institute of Mathematics USC RAS, Chernyshevskii str., 112,\\
Ufa, 450000, Russia\\
\& Bashkir State Pedagogical University, October rev. st., 3a,
\\
Ufa, 450008, Russia,\\
E-mail: borisovdi@yandex.ru\\
matem.anrb.ru \& www.bspu.ru}
\author{M.~Tautenhahn}
\address{Fakult\"at f\"ur Mathematik, Reichenhainer Str. 41,\\
Chemnitz, D-09126, Germany,\\
www.tu-chemnitz.de/\~{}mtau}
\author{I.~Veseli\'c}
\address{Fakult\"at f\"ur Mathematik, Reichenhainer Str. 41,\\
Chemnitz, D-09126, Germany,\\
www.tu-chemnitz.de/stochastik/}
\begin{abstract}
We discuss properties of $L^2$-eigenfunctions
of Schr\"odinger operators and elliptic partial differential
operators.
The focus is set on unique continuation principles
and equidistribution properties.
We review recent results and announce new ones.
\end{abstract}
\bodymatter
\section{Introduction}
In this note we present recent results in Harmonic Analysis for solutions
of (time-independent) Schr{\"o}dinger equations and other partial differential equations.
They are motivated by interest in techniques relevant for proving
localization for random Schr{\"o}dinger operators.
The mentioned Harmonic Analysis results which we present
are  a quantitative unique continuation principle and
an equidistribution property for eigenfunctions, which is scale-uniform.
These results, and variants thereof, go under various names, depending on the particular
field of mathematics: They are called observability estimate, uncertainty relation,
scale-free unique continuation principle, or local positive definiteness.
The latter term signifies that a self-adjoint operator
is (strictly) positive definite
when restricted to a relevant subspace, while it is not so on the whole
Hilbert space.
For the purpose of motivation we discuss this property in the next section.
\par
The term \emph{localization} refers to the phenomenon,
that quantum Hamiltonians describing the movement of electrons in certain
disordered media exhibit pure point spectrum in appropriately specified energy regions.
The corresponding eigenfunctions decay exponentially in space. The (time-dependent)
wavepackets describing electrons stay localized essentially
in a compact region of space for all times. Nota bene, all mentioned properties hold
\emph{almost surely}. This is natural in the context of random operators.
\par
An important partial result
for deriving localization are Wegner estimates.
These are bounds on the expected number of eigenvalues in a bounded
  energy interval of a random Schr{\"o}dinger operator restricted to a box.
\par
The localization problem has been studied for other classes of random operators
beyond those of Schr{\"o}dinger type. An example are random divergence type operators,
see e.g.\ Refs.~\citenum{FigotinK-96} and \citenum{Stollmann-98}. 
This are partial differential operators with randomness in coefficients
of higher order terms.
In paricular, the second order term is no longer the Laplacian,
but a variable coefficient operator.
In this context one is again lead to consider the above mentioned questions
of Harmonic Analysis for eigenfunctions of differential operators.
In this note we present an exposition of recently published results, and an announcement of a  quantitative unique continuation principle and
an equidistribution estimate  for eigenfunctions
for a class of elliptic operators with variable coefficients.
\subsection{Motivation: Moving and lifting of eigenvalues}
\label{ss:motivation}
Here we discuss some aspects of eigenvalue perturbation theory.
It will provide an accessible explanation why one is interested
in the results presented in Sections~\ref{sec:schroedinger} and \ref{sec:elliptic}  below
in the context of random Schr{\"o}dinger operators and elliptic differential operators, respectively.
In fact, to illustrate the main questions it will be for the moment
completely sufficient to restrict our attention to the finite dimensional situation,
i.e.\ to perturbation theory for finite symmetric matrices.
The focus will be on how
(local) positive definiteness of the perturbation relates to lifting of eigenvalues.
\par
Let $A$ and $B$ be symmetric $n \times n$ matrices, with $B\geq b>0$ positive definite.
The variational min-max principle for eigenvalues
shows that for any $k \in \{1,\dots, n\}$ and $t\geq 0$
\begin{equation}
\label{eq:positive_definite_perturbation}
\lambda_k(A+tB) \geq \lambda_k(A) + b \, t
\end{equation}
where $\lambda_k(M)$ denotes the $k$th lowest eigenvalue,
counting multiplicities,
of a symmetric matrix $M$.
Note that the dimension $n$ does not enter in the bound
\eqref{eq:positive_definite_perturbation}.
Without the positive definiteness assumption
on $B$ this universal bound will fail, most blatantly if
\begin{equation*}
A =\begin{pmatrix}
A_1 & \ 0 \\
0 & A_2
\end{pmatrix}
\quad \text{and} \quad
B =\begin{pmatrix}
\operatorname{Id} & \ 0    \\
0 & -\operatorname{Id}
\end{pmatrix} .
\end{equation*}
In this case, all eigenvalue $\lambda_k(A+tB)$ \emph{will move},
even with constant speed w.r.t.\ the variable $t$, albeit in different directions.
If $B$ is singular, some eigenvalues may not move at all.
However, for appropriate classes of symmetric matrices $A$,
and of positive semidefinite matrices $B$, one may still aim to prove
\begin{equation}
\label{eq:positive_semidefinite_perturbation}
\forall \, t\geq 0, k \in \{1,\dots, n\} \, \exists \,  \kappa >0 \text{ such that }
\lambda_k(A+tB) \geq \lambda_k(A) + \kappa  t
\end{equation}
Note however, that $\kappa $ is now not a uniform bound but
depends on
\begin{itemlist}
 \item
the class of symmetric matrices from which $A$ is chosen,
 \item
the class of semidefinite matrices from which $B$ is chosen,
 \item
the range from which the coupling $t$ is chosen, and
 \item
the range from which the index $k \in \{1,\dots, n\}$ is chosen.
\end{itemlist}
In the case of random operators or matrices one in is interested in the situation where
\begin{equation}
\label{eq:multiparameter_family}
A(\omega)
=A_0+\sum_{j \in Q} \omega_j B_j
=\Big(A_0+\sum_{j \in Q, j\neq 0} \omega_j B_j \Big)+ \omega_0 B_0
\end{equation}
is a multi-parameter pencil.
Here $Q$ is some subset of $\ZZ^d$ containing $0$.
The real variables $\omega_j$
model random coupling constants determining
the strength of the perturbation $B_j$ in each configuration
$\omega=(\omega_j)_{j \in Q}$.
Now, \eqref{eq:multiparameter_family}
already suggest to write $A(\omega)$ as
\[
A(\omega_0^\perp)+tB \quad\text{where} \quad
t=\omega_0,\ B=B_0,\ \text{and} \ \omega_0^\perp=(\omega_j)_{j \in Q, j\neq0} .
\]
This highlights that if we consider $A(\omega)$ as a function of the single
variable $t=\omega_0$, it is clearly a one-parameter family of operators,
albeit the ``unperturbed part'' $A(\omega_0^\perp)$ of
$A(\omega)=A(\omega_0^\perp)+tB$ is not a single operator,
but varying over the ensemble $(A(\omega_0^\perp))_{\omega_0^\perp}$.
To have a useful version of \eqref{eq:positive_semidefinite_perturbation}
in this situation, the constant $\kappa $ needs to have a uniform
lower bound $\inf_{A} \kappa $ where $A=A(\omega_0^\perp)$
varies over all matrices in the ensemble.
\par
In what follows we present
rigorous results of the type \eqref{eq:positive_semidefinite_perturbation},
but where $A$ and $B$ are not finite matrices, but differential and multiplication
operators.
The relevant operators have all compact resolvent, ensuring
that the entire spectrum consists of eigenvalues.
%
%
%
%
\section{Equidistribution property of Schr\"odinger eigenfunctions} \label{sec:schroedinger}
The following result is taken from Ref.~\citenum{Rojas-MolinaV-13}.
It is an equidistribution estimate for Schr{\"o}dinger eigenfunctions, which is uniform
w.r.t.\ the naturally arising length scales, and has strong implications for
the spectral theory of random Schr\"odinger operators.
\par
We fix some notation. For $L>0$ we denote by $\Lambda_L = (-L/2 , L/2)^d$ a cube in $\RR^d$.
For $\delta>0$ the open ball centered at $x\in \RR$ with radius $\delta$ is denoted by
$B(x, \delta)$. For a sequence of points $(x_j)_j$ indexed by $j \in \ZZ^d$
we denote the collection of balls  $\cup_{j \in \ZZ^d} B(x_j , \delta) $
by $S$ and its intersection with $\Lambda_L$ by $S_L$.
%
%
We will be dealing with certain subspaces of the standard second order Sobolev space
$W^{2,2}(\L_L)$ on the cube.
Let  $\Delta$ be the $d$-dimensional Laplacian.
Its restriction to the cube $\L=\L_L$ needs boundary conditions to be self-adjoint.
The domain of the Dirichlet Laplacian will be denoted by
$\cD(\Delta_{\L,0})$ and the domain of the Laplacian with periodic boundary conditions
by $\cD(\Delta_{\L,\mathrm{per}})$.
Let $V \colon \RR^d\to \RR$ be a bounded measurable
function, and $H_L = (-\Delta + V)_{\Lambda_L} $ a Schr\"odinger operator on the cube
$\Lambda_L$ with Dirichlet or periodic boundary conditions.
The corresponding domains are still
$\cD(\Delta_{\L,0})$ and $ \cD(\Delta_{\L,\mathrm{per}})$, respectively.
Note that we denote a multiplication operator by the same symbol as the corresponding function.
\par
The following theorem was proven in Ref.~\citenum{Rojas-MolinaV-13}.
\begin{theorem}[Scale-free unique continuation principle]
\label{thm:RojasVeselic}
 Let $\delta, K_{V} > 0$. Then there exists $C_{\rm sfUC} \in (0,\infty)$ such that for all
 $L \in 2\NN+1 $, all measurable   $V : \RR^d \to [-K_{V} , K_{V}]$, all real-valued $\psi \in
 \cD(\Delta_{\L,0})  \cup  \cD(\Delta_{\L,\mathrm{per}})$
 with $(-\Delta + V)\psi = 0$ almost everywhere on $\Lambda_L$, and  all sequences $(x_j)_{j \in \ZZ^d} \subset \RR^d$,
 such that  for all $j \in \ZZ^d$ the ball $B(x_j , \delta) \subset \Lambda_1 + j$, we have
\begin{equation}
\label{eq:observability}
\int_{S_L} \psi^2  \geq C_{\rm sfUC}  \int_{\Lambda_L} \psi^2 .
\end{equation}
\end{theorem}
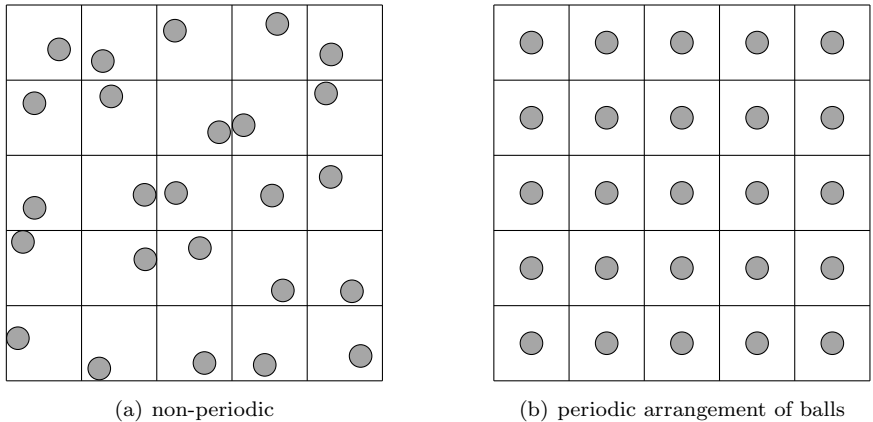
\begin{figure}[h]\centering
\subfigure[non-periodic]{
\begin{tikzpicture}
\pgfmathsetseed{{\number\pdfrandomseed}}
\foreach \x in {0.5,1.5,...,4.5}{
  \foreach \y in {0.5,1.5,...,4.5}{
    \filldraw[fill=gray!70] (\x+rand*0.35,\y+rand*0.35) circle (0.15cm);
  }
}
\foreach \y in {0,1,2,3,4,5}{
\draw (\y,0) --(\y,5);
\draw (0,\y) --(5,\y);
}
\end{tikzpicture}
}
\hspace{1cm}
\subfigure[periodic arrangement of balls]{
\begin{tikzpicture}
\foreach \x in {0.5,1.5,...,4.5}{
  \foreach \y in {0.5,1.5,...,4.5}{
    \filldraw[fill=gray!70] (\x,\y) circle (1.5mm);
  }
}
\foreach \y in {0,1,2,3,4,5}{
  \draw (\y,0) --(\y,5);
  \draw (0,\y) --(5,\y);
}
\end{tikzpicture}
}
\caption{Examples of collections of balls $S_L$ within region $\Lambda_L\subset \RR^2$.\label{fig:equidistributed}}
\end{figure}
%
%
The value of  the result is not in the \emph{existence}
of the constant $C_{\rm sfUC}$,  but in the \emph{quantitative control}
of the dependence of $C_{\rm sfUC}$ on parameters entering the
model. The very formulation of the theorem  states that $C_{\rm sfUC}$ is
independent of the position of the balls $B(x_j,\delta)$ within $\Lambda_1 +j$,
and independent of the scale $L\in2\mathbb{N} +1$.
From the estimates given in Section~2 of Ref.~\citenum{Rojas-MolinaV-13} one infers
that $C_{\rm sfUC}$ depends on the potential
$V$ only through the norm $\lVert V \rVert_\infty$ (on an exponential scale),
and it depends on the small radius $\delta>0$ polynomially, i.e.\ $C\gtrsim \delta^N$,
for some $N\in\mathbb{N}$ which depends on the dimension on $d$ and $\lVert V \rVert_\infty$.
\par
The theorem states a property of functions in the kernel of the operator.
It is easily applied to eigenfunctions corresponding to other eigenvalues since
\[
 H_L\psi=E\psi \Leftrightarrow (H_L-E)\psi=0 .
\]
As a consequence of the energy shift the constant $K_{V}$ has to be replaced
with $K_{V-E}$, which may be larger than $K_{V}$. It may always be estimated by
 $K_{V-E}\leq K_V+|E|$.
\par
There is a very natural question supported by earlier results,
which was spelled out in Ref.~\citenum{Rojas-MolinaV-13},
namely does the following generalisation of Theorem~\ref{thm:RojasVeselic} hold:
Given $\delta >0$, $K\geq0$ and $E\in\mathbb{R}$ there is a constant $C>0$ such that
for all measurable $ V\colon \mathbb{R}^d \rightarrow [-K,K] $, all $L \in 2\NN+1$,
and all sequences $(x_j)_{j\in\mathbb{Z}^d} \subset \mathbb{R}^d$
with $B(x_j,\delta)  \subset\Lambda_1 +j$ for all $j \in \ZZ^d$ we have
\begin{equation}
\label{eq:uncertainty}
 \chi_{(-\infty,E]} (H_L) \, W_L \, \chi_{(-\infty,E]} (H_L) \geq C~ \chi_{(-\infty,E]} (H_L) ,
\end{equation}
where $W_L=\chi_{S_L}$  is the indicator function of  $S_L$ and
$\chi_{I} (H_L)$ denotes the spectral projector of $H_L$ onto the interval $I$.
Here $C=C_{\delta, K, E}$ is determined by $\delta, K, E$ alone.
\par
Klein obtained a positive answer to the question for sufficiently
short subintervals of $(-\infty,E]$.
\begin{theorem}[Ref.~\citenum{Klein-13}] \label{thm:Klein-13}
Let $d \in \NN$, $E\in \RR$, $\delta\in (0,1/2]$ and $V:\RR^d \to \RR$ be measurable and bounded. There is a constant $M_d>0$ such that if we set
\[
 \gamma = \frac{1}{2} \delta^{M_d \bigl(1 + (2\lVert V \rVert_\infty + E)^{2/3}\bigr)} ,
\]
then for all energy intervals $I\subset (-\infty, E]$ with length bounded by $2\gamma$, all $L \in 2\NN+1$, $L\geq 72 \sqrt{d}$
and all sequences $(x_j)_{j\in\mathbb{Z}^d} \subset \mathbb{R}^d$
with $B(x_j,\delta)  \subset\Lambda_1 +j$ for all $j \in \ZZ^d$
\begin{equation}
 \chi_{I} (H_L) \, W_L \, \chi_{I} (H_L) \geq \gamma^2\chi_{I} (H_L) .
\end{equation}
\end{theorem}
This does not answer the above posed question question completely due to the restriction
$|I| \leq 2\gamma$.
However, the result is sufficient for many questions in spectral theory of random Schr\"odinger operators.
For a history of the questions discussed here and earlier results we refer to Ref.~\citenum{Rojas-MolinaV-13}.
\subsection{Random Schr{\"o}dinger operators}\label{ss:rSo}
Let ${\L_L}$ be a cube of side $L\in2\NN+1$, $(\Omega, \PP)$ a probability space,
$V_0 \colon {\L_L}\to \RR$ a bounded, measurable deterministic potential, $V_\omega \colon {\L_L}\to \RR$ a bounded random
potential and $H_{\omega,L}= (-\Delta + V_0+V_\omega)_{\L_L}$ a random Schr\"odinger operator on $L^2({\L_L})$ with Dirichlet or periodic boundary conditions.
We assume that the random potential is of Delone-Anderson form
\begin{equation*}
 V_\omega(x):= \sum_{j \in{\ZZ^d}} \ \omega_j u_j(x) .
\end{equation*}
The random variables $\omega_j, j\in {\ZZ^d},$ are independent with probability distributions $\mu_j$, such that
for some $m>0$ an all $j\in {\ZZ^d}$ we have $\supp \mu_j \subset [-m, m]$.
Fix $0 < \delta_- < \delta_+<\infty$ and $0 < C_- \leq C_+ <\infty$.
The sequence of measurable functions $u_j \colon \RR^d \to \RR$, $j \in {\ZZ^d}$, is
such that
\begin{align*}
\forall j \in {\ZZ^d}:
 \quad  C_- \chi_{B(z_j,\delta_-)} \leq u_j \leq C_+ \chi_{B(z_j,\delta_+)}, \
\text{and} \ B(z_j,\delta_-) \subset \L_1 + j .
 \end{align*}
\subsection{Lifting of eigenvalues}
\label{ss:lifting}
Let $\lambda_k^L(\omega)$ denote the eigenvalues of $H_{\omega,L}$
enumerated in non-decreasing order and counting multiplicities
and $\psi_k=\psi_k^L(\omega)$ the  normalised eigenvectors corresponding to $\lambda_k^L(\omega)$.
While we suppress the dependence of $\psi_k$ on $L$ and $\omega$ in the notation,
it should be kept in mind.
Then
\[
\lambda_k^L(\omega)
=
\la \psi_k, H_{\omega,L} \psi_k\ra
=
\int_{\L_L} \overline{\psi_k} ( H_{\omega,L} \psi_k ) .
\]
Define the vector $ e=(e_j)_{j\in{\ZZ^d}}$ by  $e_j=1$ for $j\in{\ZZ^d}$.
Consider the monotone shift of $V_\omega$
\[
 V_{\omega+ {t} \cdot e}
= \sum_{j \in{\ZZ^d}} (\omega_j+ {t} ) u_j
\]
and set $\tL=\tL_L= \Lambda_L \cap {\ZZ^d}$.
By first order perturbation theory we have
\[
\frac{\rm d}{{\rm d}{\tau}} \lambda_k^L(\omega+ {\tau} \cdot e) |_{\tau=t}
= \langle \psi_k, \sum_{k \in \tL}  u_j \,   \psi_k \ra.
\]
Note that the right hand side depends on $t$ implicitly through the eigenfunction $\psi_k$.
Let us fix some $E_0\in\RR$ and restrict our attention only to those eigenvalues
satisfying $\lambda_n^L(\omega) \leq E_0$.
By Theorem~\ref{thm:RojasVeselic} there exists a constant $C_{\rm sfUC}$ depending on the energy $E_0$,
$\delta_-$ and the overall supremum
\begin{equation*}
 \label{eq:Vsupremum}
\sup_{|s|\leq m} \ \sup_{|\omega_j|\leq m} \ \sup_{x\in\RR^d} 
\big|V_{0}(x) +V_\omega(x) +s \sum_{j\in Q} u_j \big|
\end{equation*}
of the potential, such that
\begin{equation*}
\sum_{k \in \tL}  \la \psi_k, u_j \,  \psi_k \ra
\geq
C_- \sum_{k \in \tL}\la \psi_k,  \chi_{B(z_k,\delta_-)}\psi_k \ra \geq
C_-\cdot C_{\rm sfUC}
=: \kappa .
\end{equation*}
Here we used that $\|\psi\|_{L^2 (\Lambda)}=1$.
(Note that the quantity $\kappa$ depends a-priori on the model parameters.)
Integrating the derivative gives
\begin{align}
 \nonumber
\lambda_k^L(\omega+ {t} \cdot e)
&=
\lambda_k^L(\omega) + \int_0^{t}  \frac{\drm \lambda_k^L(\omega+ \tau \cdot e) }{\drm \tau}|_{\tau=s} \, \drm s \\
& \geq
\lambda_k^L(\omega) + \int_0^{t}  \kappa \, \drm s  = \lambda_k^L(\omega) + t  \kappa .
\label{eq:lifting}
\end{align}
This is the lifting estimate for eigenvalues of random (Schr\"odinger) operators
alluded to in \S \ref{ss:motivation}. It should be compared with
\eqref{eq:positive_semidefinite_perturbation} there.
Indeed, due to the uniform nature of the estimate in
Theorem~\ref{thm:RojasVeselic} we have 
\begin{equation}
\label{eq:uniform_kappa}
\inf_{ L \in 2\NN+1}
\ \inf_{\omega \text{ s.t. } \forall \, k : |\omega_j|\leq m}
\ \inf_{ |{t}|\leq m}
\ \inf_{n \text{ s.t. } \lambda_n^L(\omega)\leq E_0} \kappa >0 .
\end{equation}
Thus eigenvalues lifting estimate is almost as uniform as
\eqref{eq:positive_definite_perturbation}. A parameter, with
respect to which the lifting estimate is \emph{not} uniform is the cut-off energy $E_0$.
Indeed, if we add in \eqref{eq:uniform_kappa} an infimum over $E_0>0$ on the left hand side,
it becomes zero,
unless $\sum_k\chi_{B(z_k,\delta_-)}\geq 1$ almost everywhere on $\RR^d$.
\subsection{Wegner estimates}
Here we present a Wegner estimate. Such estimates play an important role in the proof of localization
via the multiscale analysis. The latter is an induction argument over increasing
length scales. The Wegner bound is used to prove the induction step.
\par
Let $ s\colon [0,\infty) \to[0,1]$ be the global modulus of continuity of the family
$\{\mu_j\}_{j\in {\ZZ^d}}$, that is,
\begin{equation*}
\label{definition-s-mu-epsilon}
 s(\epsilon):= \sup_{j \in {\ZZ^d}}
\sup_{a \in \RR} \, \mu_j\Big(\Big[a-\frac{\epsilon}{2},a+\frac{\epsilon}{2}\Big]\Big)
\end{equation*}
The main result of Ref.~\citenum{Rojas-MolinaV-13} on the model described in the last paragraph is a
Wegner estimate which is valid for all compact energy intervals.
\begin{theorem}[Ref.~\citenum{Rojas-MolinaV-13}]
\label{t:Wegner}
Let $H_{\omega,L}$ be a random Schr\"odinger operator as in \S \ref{ss:rSo}.
Then for each $E_0\in \RR$ there exists a constant $C_W$, such that for all $E\le E_0$,
$\epsilon \le 1/3$, and all $L\in 2\NN+1$ we have
\begin{equation*}
\label{eq:WE}
\EE\{\Tr [ \chi_{[E-\epsilon,E+\epsilon]}(H_{\omega, L}) ]\}
\le C_W \ s(\epsilon) \, \lvert \ln \, \epsilon \rvert^d \ \lvert \Lambda_L \rvert .
\end{equation*}
 \end{theorem}
The Wegner constant $C_W$ depends only on $E_0$, $\|V_0\|_\infty$, $m$,  $C_-$, $C_+$, $\delta_-$,
and  $\delta_+$.
Klein\cite{Klein-13} obtains an improvement over this result based on his
above quoted Theorem~\ref{thm:Klein-13}. There are many earlier, related Wegner estimates.
For an overview we refer to Ref.~\citenum{Rojas-MolinaV-13}.
%
%
%
%
%
\subsection{Comparison of local $L^2$-norms}
An important step in the proof of Theorem~\ref{thm:RojasVeselic}
is the following result which compares $L^2$-norms of the restrictions
of a PDE-solution to two distinct subsets. In our applications
the solution will be an eigenfunction of the Schr\"odinger operator.
Various estimates of this type have been given in Refs.~\citenum{GerminetK-13}, \citenum{BourgainK-13} and \citenum{Rojas-MolinaV-13}. 
We quote here the version from the last mentioned paper.
\begin{theorem}
\label{thm:quantitative-UCP}
Let $K, R, \beta\in [0, \infty), \delta \in (0,1]$.
There exists a constant $C_{\rm qUC}=C_{\rm qUC}(d,K, R,\delta, \beta) >0$ such that,
for any $G\subset  \RR^d$ open, any $\Theta\subset G$ measurable,
satisfying the geometric conditions
\[
\operatorname{diam} \Theta +  \operatorname{dist} (0 , \Theta) \leq 2R \leq 2  \operatorname{dist} (0 , \Theta), \quad \delta < 4R, \quad B(0, 14R ) \subset G,
\]
and any measurable $V\colon G \to [-K,K]$  and real-valued $\psi\in W^{2,2}(G)$ satisfying the  differential inequality
\begin{equation*}
\label{eq:subsolution}
\lvert \Delta \psi \rvert \leq \lvert V\psi \rvert \quad \text{a.e.on }   G
\quad  \text{ as well as } \quad
\int_{G}  \lvert \psi \rvert^2
\leq
\beta \int_{\Theta} \lvert \psi \rvert^2 ,
\end{equation*}
 we have
\begin{equation}
\label{eq:aim}
\int_{B(0,\delta)}  \lvert \psi \rvert^2
\geq
 C_{\rm qUC}
\int_{\Theta} \lvert \psi\rvert^2 .
\end{equation}
\end{theorem}
\begin{figure}[ht]\centering
 \begin{tikzpicture}
  \filldraw[black!10] plot[smooth cycle] coordinates{(1,1) (1,6) (4,6) (8,7.5) (10,7) (9,2) };
  \draw plot[smooth cycle] coordinates{(1,1) (1,6) (4,6) (8,7.5) (10,7) (9,2) };
  \filldraw[black!30] (5.5,3.5) circle (0.5cm);
  \draw  (5.5,3.5) circle (0.5cm);
  \filldraw (5.5,3.5) circle (1pt);
  \draw (4.4,3.7) node {\small $B(0,\delta)$};
  \filldraw[black!30]   (7,4) rectangle (8,5);
   \draw  (7,4) rectangle (8,5);
  \draw (8.3,4.5) node {\small $\Theta$};
  \draw[latex-latex] (5.51,3.51)--(7,4.1);
  \draw[latex-latex] (5.5,3.48)--(5.8,1.095);
  \draw (6.35,4.1) node {\small $R$};
  \draw (5.3,2) node {\small $14R$};
  \draw (4,6.3) node {\small $G$};
\end{tikzpicture}
\caption{Assumptions in Theorem~\ref{thm:quantitative-UCP} on the geometric constellation of $G$, $\Theta$, and $B(0,\delta)$}
\end{figure}
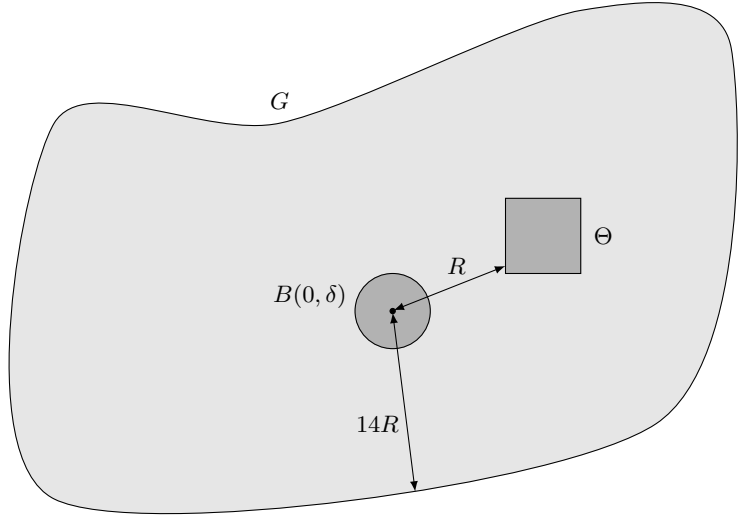
%
%
%
%
%
\section{Equidistribution property eigenfunctions of second order elliptic operators}\label{sec:elliptic}
\subsection{Notation}
Let $\cL$ be the second order partial differential operator
\[
 \cL u = -\sum_{i,j=1}^d \partial_i \left( a^{ij} \partial_j u \right)
\]
acting on functions $u$ on $\RR^d$. Here $\partial_i$ denotes the $i$th weak derivative. Moreover, we introduce the following assumption on the coefficient functions $a^{ij}$.
\begin{assumption}\label{ass:elliptic+}
Let $r,\elliptic_1 , \elliptic_2 > 0$. The operator $\cL $ satisfies $A(r,\elliptic_1 , \elliptic_2)$, if and only if $a^{ij} = a^{ji}$ for all $i,j \in \{1,\ldots , d\}$ and for almost all $x,y \in B(0,r)$ and all $\xi \in \RR^d$ we have
\begin{equation*} \label{eq:elliptic}
\elliptic_1^{-1} \lvert \xi \rvert^2 \leq \sum_{i,j=1}^d a^{ij} (x) \xi_i \xi_j \leq \elliptic_1 \lvert \xi \rvert^2 \quad\text{and}\quad \sum_{i,j=1}^d \lvert a^{ij} (x) - a^{ij} (y) \rvert \leq \elliptic_2 \lvert x-y \rvert .
\end{equation*}
\end{assumption}
\subsection{A quantitative unique continuation principle}
We first present an extension of the quantitative
continuation principle, formulated for Schr\"odinger operators
in Theorem~\ref{thm:quantitative-UCP},
to elliptic operators with variable coefficients.
\begin{theorem}[Ref.~\citenum{BorisovTV}] \label{thm:qUC-elliptic}
Let $R\in (0,\infty)$, $K_V, \beta \in [0,\infty)$ and $\delta \in (0, 4 R]$.
There is an  $\epsilon> 0$, such that if $ A(14R, 1+\epsilon, \epsilon)$ holds
then there is a constant $C_{\rm qUC} > 0$, such that for any open $G\subset \RR^d$
containing the origin and $\Theta \subset G$ measurable satisfying
 \[
 \operatorname{diam} \Theta + \dist (0 , \Theta) \leq 2R \leq 2 \dist (0 , \Theta) \quad
 \text{and} \quad  B(0,14R) \subset G,
\]
any measurable $V : G \to [-K_V , K_V]$ and real-valued $\psi \in W^{2,2} (G)$ satisfying the differential inequality
\begin{equation*} \label{eq:psi}
 \lvert \cL  \psi \rvert \leq \lvert V\psi \rvert \quad \text{a.e.\ on $G$} \quad \text{as well as} \quad \frac{\lVert \psi \rVert_G^2}{\lVert \psi \rVert_\Theta^2} \leq \beta ,
\end{equation*}
we have
\begin{equation}
 \lVert \psi \rVert_{B(x,\delta)}^2 \geq C_{\rm qUC} \lVert \psi \rVert_{\Theta}^2 .
\end{equation}
\end{theorem}
\subsection{Scale-free unique continuation principle}
We move on to discuss the equidistribution property or
scale-free unique continuation principle for eigenfunctions.
The aim is to formulate an analog of Theorem~\ref{thm:RojasVeselic}
for variable coefficient elliptic operators. As presented below,
for the moment we have solved only the situation where the second
order term is sufficiently close to the Laplacian.
\par
As before, we denote by $\L_L$ a box of side $L\in \NN$.
By $V$ we indicate a bounded measurable potential on
$\RR^d$ taking values in $[-K_V,K_V]$,
where $K_V$ is a positive constant.
We restrict the operator $\cL $ on $\L_L(0)$ and add either periodic or Dirichlet boundary conditions.
In the former case we denote such an operator by $\cL _{L,0}$,
and its domain $\Dom(\cL _{L,0})$ is the subspace of $W^{2,2}(\L_L)$
consisting of functions vanishing on $\partial \L_L$.
The notation for the operator with periodic boundary condition
is $\cL _{L,\mathrm{per}}$ and its domains $\Dom(\cL _{L,\mathrm{per}})$
consists of the functions in $W^{2,2}(\L_L)$ satisfying periodic boundary conditions.
\begin{assumption}\label{ass:periodicCoefficients}
For each pair $i,j$ the function $a^{ij}\colon \RR^d \to \RR$ is $\ZZ^d$-periodic.
\par
Assume that in the case of operator $\cL _{L,0}$ its coefficients $a^{ij}$, $i\not= j$  vanish on the sides of box $\L_L$, while the coefficients $a^{ii}$ satisfy periodic boundary conditions on the sides of box $\L_L$. In the case of operator $\cL _{L,\mathrm{per}}$ suppose that all its coefficients satisfy periodic boundary conditions on the sides of box $\L_L$.
\end{assumption}
\begin{theorem}\label{thm:equidistribution-elliptic}
Fix $K_V\in [0,+\infty)$, $\delta\in(0,1]$. Assume $A(\sqrt{d},1+\epsilon,\epsilon)$ with $\epsilon>0$ as in
Theorem \ref{thm:qUC-elliptic} . Assume \ref{ass:periodicCoefficients}.
\par
Then there exists a constant $C_{sfUC}>0$ such that for any $L\in 2\NN+1$, any sequence
\begin{equation*}\label{d1.1}
Z:=\{z_k\}_{k\in\ZZ^d} \ \text{ in }\ \RR^d
\quad \text{such that} \ B(z_k,\delta)\subset \L_1(k)  \text{ for each } k\in\ZZ^d,
\end{equation*}
any measurable $V: \L_L\mapsto [-K_V,K_V]$ and any real-valued  $\psi\in\Dom(\cL _{L,0})$, respectively $\psi\in \Dom(\cL _{L,\mathrm{per}})$ satisfying 
\begin{equation*}\label{d1.2}
|\cL\psi|\leqslant |V\psi|\quad \text{a.e.}\quad \L_L
\end{equation*}
we have
\begin{equation}\label{d1.3}
\int\limits_{S_L} |\psi(x)|^2 dx=\sum\limits_{k\in Q_L} \|\psi\|_{L_2(B(z_k,\delta))}^2\geqslant C_{sfUC} \|\psi\|_{L_2(\L_L)}^2,
\end{equation}
where $S_L:=S\cap\L_L=\cup_{k\in Q_L} B(z_k,\delta)$, $Q_L=\L_L\cap \ZZ^d$, and $S:=\cup_{k\in \ZZ^d} B(z_k,\delta)$.
\end{theorem}
As a \emph{Corollary} we obtain immediately an eigenvalue lifting estimate
analogous to \eqref{eq:lifting}, where $\kappa$ is again uniform w.r.t.\ many parameters,
as spelled out in subsection \ref{ss:lifting} explicitly.
\par
The proof of Theorem~\ref{thm:equidistribution-elliptic}
is based on the strategy implemented in Ref.~\citenum{Rojas-MolinaV-13}.
First one uses the conditions on the coefficients $a^{ij}$ described in Assumption
\ref{ass:periodicCoefficients} to extend $\psi$ as well as the differential expression
$\cL$ to the whole of $\RR^d$ while keeping the $W^{2,2}$-regularity and
the differential inequality originally satisfied by $\psi$.
Then one uses the comparison Theorem~\ref{thm:qUC-elliptic} for local $L^2$-norms.
Note that now the condition concerning the minimal distance to the boundary of $G$ plays no role,
since $\psi$ has been extended to the whole of $\RR^d$. From this point the
combinatorial and geometric arguments of Ref~\citenum{Rojas-MolinaV-13} take over.
In fact, one can prove a abstract meta-theorem:
Once the comparison of local $L^2$-norms of $\psi$ holds up to the boundary,
an equidistribution property for $\psi$ follows.
Interestingly, such an argument no longer uses the fact
that $\psi$ is a solution of an differential equation or inequality.
\section*{Acknowledgments}
D.B. was partially supported by RFBR,
the grant of the President of Russia for young scientists - doctors of science (MD-183.2014.1),
and the fellowship of Dynasty foundation for young mathematicians.
\par
M.T. and I.V. have been partially supported by the DAAD and the Croatian Ministry of Science, Education and Sports
through the PPP-grant ``Scale-uniform controllability of partial differential equations''.
M.T. and I.V. have been partially supported by the DFG.
\end{document}